\documentclass[11pt]{amsart}
\usepackage{amsfonts}
\usepackage{hyperref}
\usepackage[utf8]{inputenc}
\usepackage[T1]{fontenc}
\usepackage[dvips]{graphicx}
\usepackage{amssymb}
\usepackage{amsmath}
\usepackage{dsfont}
\usepackage{color}
\usepackage{latexsym}
\usepackage{bbm}
\usepackage{color}
\usepackage{amsthm}
\usepackage{multicol}
\usepackage[top   = 2.75cm,
bottom = 2.50cm,
left   = 2.50cm,
right  = 2.00cm]{geometry}

\bibstyle{plain}
\theoremstyle{plain}
\newtheorem{theorem}{Theorem}

\newtheorem{problem}{Problem}
\newtheorem{corollary}[theorem]{Corollary}

\newtheorem*{GNCT}{Guthrie-Nymann Classification Theorem}
\theoremstyle{remark}

\begin{document}

\title{Achievement sets \--- current results and open problems
}
\author[Sz. G\l\c{a}b and F. Prus-Wi\'{s}niowski]{ Szymon G\l\c{a}b and Franciszek Prus-Wi\'{s}niowski}

\newcommand{\eacr}{\newline\indent}

\address{\llap{*\,}Szymon G\l\c{a}b\eacr
Institute of Mathematics\eacr
\L\'{o}d\'{z} University of Technology\eacr
ul. W\'{o}lcza\'{n}ska 215\eacr
PL-93-005 \L\'{o}d\'{z}\eacr
Poland\acr ORCID 0000-0001-9026-8235}
\email{szymon.glab@p.lodz.pl}

\address{\llap{*\,}Franciszek Prus-Wi\'{s}niowski\eacr
Instytut Matematyki\eacr
Uniwersytet Szczeci\'{n}ski\eacr
ul. Wielkopolska 15\eacr
PL-70-453 Szczecin\eacr
Poland\acr ORCID 0000-0002-0275-6122}
\email{franciszek.prus-wisniowski@usz.edu.pl}

\subjclass[2020]{Primary: 40A05; Secondary: 11B05, 28A75}
\keywords{ Achievement set, Cantor set, Cantorval, center of distances}

\date{\today}

\newcommand{\acr}{\newline\indent}

\begin{abstract}
We survey recent developments in the theory of achievement sets and present a substantial collection of open problems.
\end{abstract}

\maketitle

One of the aims of this nonstandard article is to highlight the most interesting recent results in achievement set theory. Its primary focus, however, is on open problems in the area. Almost all of the problems presented here have been known for some time and have resisted several attempts at resolution.

The idea of writing the article  was inspired by two papers of the late Dan Waterman \cite{W1}, \cite{W2}. Those works stimulated renewed interest in his main area of research—bounded generalized variation—and led to a number of further developments. We cherish the memory of this fine mathematician and, above all, of a man of great wisdom and kindness.

\section{The Kakeya conditions}

The first mathematician to investigate sets of subsums of a convergent series in a systematic way was Soichi Kakeya, more than a century ago \cite{Kakeya,Kakeya2}. He was interested in the topological nature of such sets, and one of his pioneering results states that the \textsl{set of subsums} (equivalently, the \textsl{achievement set})
\[
E = E(x_n) := \Bigl\{ y \in \mathbb{R} : \exists\, A \subset \mathbb{N} \text{ such that } y = \sum_{n \in A} x_n \Bigr\}
\]
of an absolutely convergent real series $\sum x_n$ is always compact and perfect. In fact, the above definition makes sense for any real sequence $(x_n)$ (see \cite{Jones}), and over time the term \textsl{achievement set} has prevailed over the older name \textsl{set of subsums}. One of the fundamental results in the theory of achievement sets is the celebrated Guthrie–Nymann Classification Theorem \cite[Thm.~1]{GN88}.

\begin{GNCT}\label{t6p2}
The set $E$ of all subsums of an absolutely convergent series is always of exactly one of the following four types:
\begin{itemize}
\item[(i)] a finite set;
\item[(ii)] a multi-interval set;
\item[(iii)] a Cantor set;
\item[(iv)] a Cantorval.
\end{itemize}
\end{GNCT}

A \textsl{multi-interval set} is the union of finitely many bounded closed intervals. A \textsl{Cantor set} is a set homeomorphic to the classical ternary Cantor set. \textsl{Cantorvals} are the most intricate achievement sets: they are compact sets $A \subset \mathbb{R}$ such that the endpoints of every $A$-gap\footnote{$A$-gaps are the connected components of $\mathbb{R}\setminus A$.} are accumulation points both of other $A$-gaps and of $A$-intervals\footnote{$A$-intervals are the nontrivial connected components of $A$.}. Cantorvals were introduced in \cite{MO} under the name $M$-Cantorvals. Mendes and Oliveira proved that all Cantorvals are mutually homeomorphic \cite[Appendix]{MO} (see also \cite[Cor.~12]{BKP}). Clearly, all Cantor sets are homeomorphic, while multi-interval sets are homeomorphic if and only if they have the same number of components; the same holds for finite sets. Hence, the notion of \textsl{topological type} of an achievement set is well justified.

A particularly natural example of a Cantorval was given by Guthrie and Nymann \cite[p.~326]{GN88}. We refer to it as the \textsl{model Cantorval}: it is obtained by taking the ternary Cantor set $C$ and adjoining all intervals removed at odd-numbered steps of the standard construction of $C$.

Naturally, the basic question is how one can recognize the topological type of $E(x_n)$ by examining the terms of the series $\sum x_n$. This question was already investigated by Kakeya, who studied the relationship between the terms $x_n$ and the remainders $r_n=\sum_{i>n}x_i$. Kakeya considered only convergent series with positive terms, which does not affect the generality of the discussion, since
\[
E(x_n)=E(|x_n|)+\sum_{n:\,x_n<0}x_n
\]
(see \cite{Ho41}, \cite{BFPW1}). The set of indices
\[
K(x_n):=\{n\in\mathbb N:\ x_n>r_n\}
\]
is called the set of \textsl{Kakeya conditions}. Its complement,
\[
K^c(x_n)=\{n\in\mathbb N:\ x_n\le r_n\},
\]
is called the set of \textsl{reversed Kakeya conditions}. Every subset of $\mathbb N$ belongs to exactly one of the following three simple categories: finite sets, sets with finite complement, and infinite sets with infinite complement. The first two of these categories already appeared in Kakeya’s early results on the relationship between the terms of a series and the topological nature of the corresponding achievement set.

\begin{theorem}
\label{t1}
    $E(x_n)$ is a multi-interval set if and only if $K(x_n)$ is finite.
\end{theorem}

\begin{theorem}
\label{t2}
    If $K(a_n)$ has finite complement, then $E(x_n)$ is a Cantor set.
\end{theorem}

This implication is not reversible, as demonstrated by the series $\sum b_n$ with $b_{2k} = b_{2k-1} = \frac{1}{4^k}$ for all $k$. The easiest way to see that $E(b_n)$ is a Cantor set is by applying Theorem~16 from \cite{BFPW2}. By the Guthrie--Nymann Classification Theorem, the above results yield the following simple corollary.

\begin{corollary}
\label{c3}
If $E(a_n)$ is a Cantorval, then $\mathrm{card}\,K(a_n) = \mathrm{card}\,K^c(a_n) = \infty$.
\end{corollary}

These observations uniquely determine the topological type of $E(a_n)$ when $K(a_n)$ is either finite or cofinite. Kakeya believed that $\mathrm{card}\,K(a_n) = \infty$ implies that $E(a_n)$ is a Cantor set. In \cite{Kakeya}, he candidly wrote:
\textit{``That the relation $a_n \le r_n$ fails only for an infinite number of values of $n$ seems to be the necessary and sufficient condition that the set $E(a_n)$ should be nowhere dense; but I have no proof of it.''}
This statement is known as the Kakeya 1914 hypothesis and is overshadowed by the better-known Kakeya 1917 hypothesis.

The Guthrie--Nymann Cantorval $E(c_n)$, where $c_{2n-1} = \frac{3}{4^n}$ and $c_{2n} = \frac{2}{4^n}$, provides a counterexample to this hypothesis, since $K(c_n) = 2\mathbb{N}$ (see \cite{GN88}, \cite{BPW}). This example, together with $E(b_n)$, shows that the Kakeya conditions do not uniquely determine the topological type of an achievement set, as $K(c_n) = K(b_n)$. Nevertheless, there is some truth to Kakeya’s conjecture: for every infinite set $K \subset \mathbb{N}$ with infinite complement, there exists an absolutely convergent series whose set of Kakeya conditions is $K$ and whose achievement set is a Cantor set. This nontrivial result was proved by J.~Marchwicki and P.~Miska in \cite[Thm.~1.3]{MM21}. A simpler proof appears in \cite{MPWP}, though at the cost of losing uniqueness of subsums \footnote{Uniqueness of subsums means that every subsum can be obtained by one subseries only, that is, that the function $2^\mathbb N\to E(x)n):\,A\mapsto\sum_{n\in A}x_n$ is bijective.}.

\begin{theorem}
\label{t4}
For any $K \subset \mathbb{N}$ with $\mathrm{card}\,K = \mathrm{card}\,K^c = \infty$, there exists a series $\sum a_n$ such that $K(a_n) = K$ and $E(a_n)$ is a Cantor set.
\end{theorem}

They also investigated whether the Cantor set in Theorem~\ref{t4} could be replaced by a Cantorval \cite[Problem~5.1]{MM21}. This is by no means an easy question. In all known examples of achievable Cantorvals, at least half of the Kakeya conditions (in terms of asymptotic density) are reversed. For instance, the earliest known example of an achievable Cantorval is given by the Weinstein--Shapiro series $\sum d_n$, where
\[
d_n = \frac{3}{10} \cdot \frac{9 - m}{10^k},
\]
and $(k,m)$ is the unique pair in $\mathbb{N} \times \{1,2,3,4,5\}$ such that $n = 5(k-1) + m$ \cite{WS}. Since
\[
K(d_n) = \{n \in \mathbb{N} : n \equiv 0 \ (\mathrm{mod}\ 5)\},
\]
we obtain $d(K^c(d_n)) = \frac{4}{5}$.

For the Guthrie--Nymann series $\sum c_n$, we have $d(K^c(c_n)) = \frac{1}{2}$.

Motivated by these examples, Marchwicki and Miska investigated whether a Cantorval could be obtained from a sequence with a smaller proportion of reversed Kakeya conditions. They formulated the problem in terms of lower and upper asymptotic densities ($\underline{d}$ and $\overline{d}$, respectively) and proved the following result \cite[Thm.~3.10]{MM21}.

\begin{theorem}
For any $0 < \alpha \le \beta \le 1$, there exists a series $\sum a_n$ such that $E(a_n)$ is a Cantorval and
\[
\underline{d}\,K^c(a_n) = \alpha, \qquad \overline{d}\,K^c(a_n) = \beta.
\]
\end{theorem}

Moreover, they explicitly asked whether it is possible to construct a series $\sum a_n$ such that $E(a_n)$ is a Cantorval and $d(K^c(a_n)) = 0$ \cite[Problem~5.2]{MM21}. A positive answer was given last year in \cite{PWP}. This, however, was not the final development concerning the significance of Kakeya conditions. In his invited lecture \emph{A new sufficient condition for an achievement set to be a Cantorval}, delivered during the Seventh Workshop on Postmodern Real Analysis (Będlewo Conference Center, Poland, November~3--6,~2025), Piotr Nowakowski presented an ingenious construction of a series $\sum x_n$ such that $E(x_n)$ is a Cantorval and $K(x_n)=S$ for any prescribed infinite subset $S \subset \mathbb{N}$ with infinite complement. This result answers Problem~5.1 from \cite{MM21} in the affirmative and completes the overall picture of the relationship between Kakeya conditions and topological type.

In summary, the only definitive conclusions obtainable via Kakeya conditions were already discovered by Kakeya himself. Recent results show that if $K(x_n)$ is infinite with infinite complement, then $E(x_n)$ may be either a Cantor set or a Cantorval. Consequently, Kakeya conditions do not provide more effective tools for recognizing the topological type of an achievement set than the classical Kakeya theorems, namely Theorems~\ref{t1} and~\ref{t2}. Therefore, the fundamental problem of the theory of one-dimensional achievement sets essentially remains open.

\begin{problem}
Find characterizations of the topological types of achievement sets in terms of the generating sequence $(x_n)$.
\end{problem}

Further evidence that Kakeya conditions are not a universal tool for recognizing topological type comes from a counterexample constructed by Moroz \cite{M24} to a conjecture of Pratsiovytyi and Karvatskyi \cite[Conjecture~4.1]{PK23}. Moroz's example shows that even the following three assumptions:
\begin{enumerate}
\item $a_n \le c_n \le b_n$ for all $n$;
\item $K(a_n) = K(b_n) = K(c_n)$;
\item $E(a_n)$ and $E(b_n)$ have the same topological type,
\end{enumerate}
do not imply that $E(c_n)$ has the same topological type as the other two achievement sets.

The characterization of finite achievement sets is trivial, and the characterization of multi-interval sets was already obtained by Kakeya (Theorem~\ref{t1}). The cases of Cantor sets and Cantorvals, however, remain unresolved and motivate the study of special families of series or the search for new sufficient conditions for particular topological types. The difficulty of determining the topological type is well illustrated by the following concrete problems. The first originates from a paper by Jones \cite{Jones}, which played an important role in the revival of achievement set theory about fifteen years ago.

\begin{problem}[The First Jones Problem]
What is the topological type of $E(x_n)$ for $x_n := \frac{1}{2^n} + \frac{(-1)^n}{3^n}$? (see \cite[p.~518]{Jones})
\end{problem}

This series is only a minimal modification of the elementary geometric series $\sum \frac{1}{2^n}$, whose achievement set is the entire interval $[0,1]$. Nevertheless, after fourteen years, the topological type of $E(x_n)$ remains unknown. Kakeya's theorems are not applicable here, since $K(x_n) = 2\mathbb{N}$. A similar problem was posed more recently in \cite{PR25}.

\begin{problem}
What is the topological type of $E(x_n)$ for $x_{2n-1} := \frac{1}{2^{2n-1}}$ and $x_{2n} := \frac{1}{2^{2n}+1}$?
\end{problem}

More generally, the restricted goal is to find criteria determining the topological type of achievement sets for series with $K(x_n) = 2\mathbb{N}$ or $K(x_n) = 2\mathbb{N}-1$. The two problems above are closely related to the Second Jones Problem.

\begin{problem}[The Second Jones Problem]
Does there exist a sequence $(x_n)$ such that $\lim_{n \to \infty} \frac{x_{n+1}}{x_n}$ exists and $E(x_n)$ is a Cantorval? \cite[p.~519]{Jones}
\end{problem}

Another illustration of the limitations of known sufficient conditions for Cantorvals arises from the problem of thinning out slowly convergent series. A \textsl{slowly convergent} (or \textsl{interval-filling}) sequence $(x_n)$ is defined by the condition $K^c(x_n) = \mathbb{N}$. It is easy to see that every such sequence contains a \textsl{fast convergent} subsequence $(x_{n_k})$, that is, one satisfying $K(x_{n_k}) = \mathbb{N}$. By Theorem~\ref{t2}, the corresponding achievement set $E(x_{n_k})$ is then a Cantor set. On the other hand, it is known that if $x_n = r_n$ for almost all indices $n$, then $E(x_{n_k})$ is always either a multi-interval set or a Cantor set. Consequently, thinning out $(x_n)$ to obtain a subsequence yielding a Cantorval is impossible in this case. This leads to another open problem.

\begin{problem}
Let $r_n > x_n$ for all $n$. Does there exist a subsequence $(x_{n_k})$ such that $E(x_{n_k})$ is a Cantorval?
\end{problem}

Even if we restric our attention to the very friendly family of geometric sequences, the answer to the question is unknown.
\begin{problem}
Does every geometric sequence $(q^n)_{n\in\mathbb N}$ with $q\in(\frac12,1)$ contain a subsequence $(q^{n_k})_{k\in\mathbb N}$ such that $E(q^{n_k})$ is a Cantorval?
\end{problem}

\section{Multigeometric sequences}
The most important  family of sequences investigated from the point of view of the nature of their achievement sets consists of \textbf{multigeometric sequences}. They were mentioned in both unpublished versions of the note \cite{Nit15} (see \cite[p.~26]{Nit11} and \cite[footnote on p.~8]{Nit13}). Later they were defined and studied in \cite{BFS} where  the first general sufficient condition was established
 \cite[Thm.~2.1]{BFS}.

Multigeometric sequences (or series) form a class that is relatively easy to analyze from the viewpoint of the topological classification of achievement sets. Informally, a \textsl{multigeometric sequence} is a monotone mixture of finitely many convergent geometric sequences sharing the same ratio.

Formally, let $m \in \mathbb{N}$, $q \in (0,1)$, and let $k_1 \ge k_2 \ge \ldots \ge k_m > 0$. A multigeometric sequence is defined by
\[
a_n = k_i q^j,
\]
where $(j,i)$ is the unique pair with $j \in \mathbb{N}$ and $i \in \{1,\ldots,m\}$ such that $n = (j-1)m + i$. Multigeometric sequences are denoted by
\[
(k_1,\ldots,k_m;\,q),
\]
and their achievement sets by $E(k_1,\ldots,k_m;\,q)$. Note that $E(k_1,\ldots,k_m;\,q)$ is the algebraic sum
\[
k_1E(1;q)+\dots+k_mE(1;q)
\]
of $m$ scaled copies of $E(1;q)$. For $q<\tfrac12$, the set $E(1;q)$ is a Cantor set.

Historically, the first counterexamples to Kakeya’s 1914 hypothesis were constructed as achievement sets of multigeometric sequences (see \cite{WS} for an announcement without proof and \cite{F} for a proof). Later, Jones described an infinite family of multigeometric sequences whose achievement sets are Cantorvals \cite[p.~515]{Jones}. In \cite{BFS}, the first general sufficient condition for a multigeometric series to yield a Cantorval was established \cite[Thm.~2.1]{BFS}. An important example is the Guthrie–Nymann Cantorval $E(3,2;\tfrac14)$ \cite{GN88}, \cite{BPW}, which is also the achievement set of a multigeometric sequence.

A significant advance in understanding the dependence of the topological type of achievement sets of multigeometric sequences on the ratio $q$ was made in \cite{BBFS}. The authors analyzed how the topological type varies with $q$ for fixed initial coefficients $k_1q, k_2q, \ldots, k_mq$ (see \cite[pp.~1025--1026]{BBFS}). Let $E:=E(k_1,\ldots,k_m;\,q)$. There exist constants $d_C, d_{CI}, d_I$, and $d_{PM}$, depending on $k_1,\dots,k_m$, such that
\begin{itemize}
    \item If $q<d_C$, then $E$ is a Cantor set.
    \item $q\ge d_I$ if and only if $E$ is an interval; moreover, under mild assumptions, for $q<d_I$, the set $E$ is not a finite union of intervals.
    \item If $q\ge d_{CI}$, then $E$ contains a compact interval.
    \item For almost all $q\in[d_C,d_{PM}]$ (in the sense of Lebesgue measure), the set $E$ has positive measure.
    \item There exists a decreasing sequence $(q_n)$ converging to $d_C$ such that $E$ is a Cantor set of measure zero.
\end{itemize}
Always $d_C\le d_{CI}\le d_I$ and $d_C<d_{PM}$. If $d_{CI}\le q<d_I$, then $E$ is a Cantorval. Some generalizations of these conditions to a wider class of series were presented in \cite{Kar24}. Note that for $E(1,1;\,q)$ we have $d_C=d_I$; for $E(3,2;\,q)$ we have $d_C<d_{CI}=d_I$; and for $E(4,3,2;\,q)$ we have $d_C<d_{CI}<d_I$. These examples show that for some choices of $k_1,\dots,k_m$ the set $E$ can only be an interval or a Cantor set, for some the above criteria do not decide whether $E$ is a Cantorval, and for others they do.

The most recent progress in this subarea was presented in \cite{GK25}. The main theorem of that paper \cite[Thm.~8]{GK25} provides several necessary and sufficient conditions for $E(k_1,\dots,k_m;\,q)$ (in fact, in a more general setting) to have nonempty interior when $q=d_C$. This important result was obtained using tools from ergodic theory and tiling theory. These methods, new to the study of achievement sets, open promising directions for determining topological types.

One of the best-known subfamilies of multigeometric series consists of the \textbf{Ferens series}, defined and studied in \cite{BP} and \cite{BGM18}.

\begin{theorem}
\label{Ferensser}
Let $(x_n)=(m+k-1,m+k-2,\ldots,m;\,q)$ be a Ferens sequence. Then
\begin{itemize}
\item[(i)] $E(x_n)$ is a single interval if and only if $q \ge \frac{m}{s+m}=d_I$;
\item[(ii)] If $d_{CI}=\frac{1}{s-2m+1}\le q<\frac{m}{s+m}=d_I$, then $E(x_n)$ is a Cantorval;
\item[(iii)] If $0<q<\frac{1}{s-2m+3}=d_C$, then $E(x_n)$ is a Cantor set;
\item[(iv)] There exists a strictly decreasing sequence $(q_n)$ converging to $\frac{1}{s-2m+3}$ such that $E(m+k-1,\ldots,m;q_n)$ is a Cantor set of measure zero;
\item[(v)] $E(x_n)$ has positive measure for almost all $q\in\bigl(\frac{1}{s-2m+3},\frac{1}{s-2m+1}\bigr)$, in the sense of Lebesgue measure.
\end{itemize}
\end{theorem}

The next result shows that for $q=d_C$ two distinct behaviors are possible.

\begin{theorem}
Let $E:=E(m+k-1,\ldots,m;q)$ be a Ferens achievement set with $q=\frac{1}{s-2m+3}=d_C$. Then
\begin{itemize}
    \item $E$ is a Cantor set if $m\ge 3$;
    \item $E$ is a Cantorval if $m=2$.
\end{itemize}
\end{theorem}

Apart from the interval $\bigl(\frac{1}{s-2m+3},\frac{1}{s-2m+1}\bigr)$, Ferens achievement sets are completely classified. Consequently, the remaining open questions focus on this mysterious range of parameters.

\begin{problem}
Let $(x_n)$ be a Ferens sequence. Is it true that there are only countably many values of $q$ in $\bigl(\frac{1}{s-2m+3},\frac{1}{s-2m+1}\bigr)$ for which $E(x_n)$ has measure zero?
\end{problem}

It can be shown that every Ferens achievement set is an affine self-similar set. A long-standing open problem asks whether there exists an affine self-similar set of positive Lebesgue measure and empty interior in Euclidean spaces. This question was answered affirmatively in \cite{CJPPS}, where such an example was constructed in the plane. Whether this phenomenon can occur on the real line remains unknown. A negative answer to the following problem would yield a positive answer to that question; therefore, we believe the answer is positive.

\begin{problem}
Let $(x_n)$ be a Ferens sequence. Is it true that for $q\in\bigl(\frac{1}{s-2m+3},\frac{1}{s-2m+1}\bigr)$ the set $E(x_n)$ satisfies the \textsl{Palis dichotomy}, that is, it either has measure zero or contains an interval?
\end{problem}

For a multigeometric series $(k_1,\dots,k_m;\tfrac{1}{n})$, we denote by $\Sigma$ the set of all subsums of $\{k_1,\dots,k_m\}$, that is,
\[
\Sigma=\Bigl\{\sum_{i=1}^m \varepsilon_i k_i:\ \varepsilon_i\in\{0,1\}\Bigr\}.
\]
The following result is essentially due to Nitecki \cite{Nit13, Nit15}.

\begin{theorem}\label{NiteckiThm}
Let $(k_1,\dots,k_m;\tfrac{1}{n})$ be a multigeometric series such that $n=\lvert\Sigma\rvert$. Let $\{\sigma_0<\sigma_1<\dots<\sigma_{n-1}\}$ be an enumeration of $\Sigma$. For each $i$, let $t_i\in\{0,1,\dots,n-1\}$ satisfy $\sigma_i\equiv t_i\pmod n$, and assume that $\{t_0,t_1,\dots,t_{n-1}\}=\{0,1,\dots,n-1\}$. Then
\begin{itemize}
    \item[(i)] If $\sigma_i=i\sigma_1$, then $E(k_1,\dots,k_m;\tfrac{1}{n})$ is the interval $\left[0,(k_1+\dots+k_m)\frac{n}{n-1}\right]$;
    \item[(ii)] otherwise, $E(k_1,\dots,k_m;\tfrac{1}{n})$ is a Cantorval.
\end{itemize}
\end{theorem}

Let $(k_1,\ldots,k_m;1/n)$ be a multigeometric series such that $k_1$ is divisible by $n$, say $k_1=rn$. Then
\[
E(k_1,k_2,\ldots,k_m;1/n)=\{0,r\}+E(r,k_2,\ldots,k_m';1/n).
\]
Both sets have the same topological type. Hence, in determining the type of $E(k_1,\ldots,k_m;1/n)$, we may assume that none of $k_1,\ldots,k_m$ is divisible by $n$.

As an illustration, consider the sequence $(1,8;\tfrac{1}{4})$. In this case, $\Sigma=\{0,1,8,9\}$. Observe that
\[
(1,8;\tfrac{1}{4})=\left(\tfrac14,\tfrac{8}{4},\tfrac{1}{16},\tfrac{8}{16},\tfrac{1}{64},\tfrac{8}{64},\ldots\right)
=\left(\tfrac14,2,\tfrac{1}{16},\tfrac12,\tfrac{1}{64},\tfrac18,\ldots\right).
\]
After rearranging the terms, we obtain $
\left(2,\tfrac12,\tfrac14,\tfrac18,\tfrac{1}{16},\ldots\right)$.
Therefore, $E(1,8;\tfrac14)=[0,1]\cup[2,3]$.

This example shows that the Nitecki Theorem cannot be reversed. However, note that in this case $8$ is divisible by $4$. We are now ready to formulate the following conjecture.

\begin{problem}
Let $(k_1,\ldots,k_m;\tfrac{1}{n})$ be a multigeometric series such that $k_i\in\mathbb{N}$, $\Sigma=\{\sigma_0<\cdots<\sigma_{n-1}\}$ has $n$ elements, none of $k_1,\ldots,k_m$ is divisible by $n$, and $E(k_1,\ldots,k_m;\tfrac{1}{n})$ has positive Lebesgue measure. Is it true that the set of remainders $t_i<n$ with $\sigma_i\equiv t_i\pmod n$ equals $\{0,1,\ldots,n-1\}$?
\end{problem}

\textbf{Remark.} The answer is positive when $\lvert\Sigma\rvert$ is a prime number; see \cite{GK25}.

\section{Algebraic sums and decompositions of achievement sets}

A number of new results in this subarea have been established in \cite{MNP23}. Among them, the concept of $\epsilon$-\textsl{closed} (or $\epsilon$-\textsl{tight}) subsets was introduced to provide a tool for a complete characterization of sequences whose achievement sets contain an interval \cite[Prop. 28]{MNP23}. Although this characterization can be difficult to apply in many cases, it has been successfully used for several sequences, including non-multigeometric generalized Ferens sequences \cite[Thm. 3.1]{MNP23}. The most notable result in the paper is the following \cite[Thm. 4.16]{MNP23}.

\begin{theorem}
Given any $p \ge m \ge 2$ with $p, m \in \mathbb{N} \cup \{\infty\}$, there exists an achievable Cantor set $C$ such that $C_k := \underbrace{C + \cdots + C}_{\text{$k$ times}}$ satisfies
    \begin{itemize}
        \item $C_k$ is a Cantor set for every $k < m$;
        \item $C_k$ is a Cantorval for every $k$ with $m \le k < p$;
        \item $C_k$ is an interval for every $k \ge p$.
    \end{itemize}
\end{theorem}

In particular, if $C$ is an achievable Cantor set with $m=2$ and $p=\infty$, then $C_2$ is a Cantorval, and the algebraic sum of any finite number of copies of $C$ remains a Cantorval \cite[Thm. 4.14]{MNP23}. Additionally, \cite[Prop. 4.4 and 4.5]{MNP23} provide the following observation.

\begin{theorem}
For any achievable Cantorval $D$, there exist two achievable Cantor sets $C$ and $C'$ such that $D + C$ is a Cantorval and $D + C'$ is an interval.
\end{theorem}

Several similar problems remain open.

\begin{problem}
    Does there exist, for any achievable Cantor set $C$, achievable Cantor sets $C'$ and $C''$ such that $C + C'$ is a Cantor set and $C + C''$ is a Cantorval?
\end{problem}

\begin{problem}
    Does there exist, for any achievable Cantor set $C$, achievable Cantorvals $D'$ and $D''$ such that $C + D'$ is a Cantorval and $C + D''$ is an interval?
\end{problem}

\begin{problem}
    Does there exist, for any achievable Cantorval $D$, achievable Cantorvals $D'$ and $D''$ such that $D + D'$ is a Cantorval and $D + D''$ is an interval?
\end{problem}

One of the simplest subfamilies of Cantor sets consists of central Cantor sets. A \textsl{central Cantor set} with fundamental interval $[0, y_0]$ is constructed as follows. At the initial step, we remove an open interval $(y_1, x_1)$ from $[0, y_0]$, analogous to the standard Cantor ternary construction. The only requirement is $0 < x_1 - y_1 < y_0$. Denote the remaining two closed intervals by $I_1$, each of length $y_1$. In the second step, we remove a middle open interval of equal length from each component of $I_1$, ensuring it is smaller than $y_1$, and denote the leftmost removed interval by $(y_2, x_2)$. Continuing inductively, at the $n$-th step, from each of the $2^{n-1}$ components of $I_{n-1}$, we remove a middle open interval of length $x_n - y_n$ with $0 < x_n - y_n < y_{n-1}$. The set $\bigcap_{n\in\mathbb{N}} I_n$ obtained in this manner is a \textsl{central Cantor set}. All removed intervals are centrally aligned with the intervals from which they are removed. We will, unless specified otherwise, consider central Cantor sets with fundamental interval $[0,1]$.

Each central Cantor set with fundamental interval $[0,1]$ is the achievement set of the series $\sum x_n$, where $x_n$ are the right endpoints of the leftmost intervals removed at step $n$. The series is fast convergent ($K(x_n) = \mathbb{N}$), and its remainders $r_n$ correspond to the left endpoints of the removed intervals, i.e., $r_n = y_n$ for all $n$. Conversely, any fast convergent positive series of sum $1$ has an achievement set that is a central Cantor set with fundamental interval $[0,1]$.

There is a natural bijection $A$ between the set of all fast convergent positive series of sum $1$ and sequences $(a_n) \in (0,1)^\mathbb{N}$, defined by $A(x_n) = (a_n)$ with $a_n := \frac{x_n - r_n}{r_{n-1}}$. Its inverse is $A^{-1}(a_n) = (x_n)$ with $x_n := \frac{1 + a_n}{2} \prod_{i=1}^{n-1} \frac{1 - a_i}{2}$. The $a_n$ are called \textsl{ratios of dissection}, and $C_{(a_n)}$ denotes the central Cantor set corresponding to $(a_n)$. Central Cantor sets are symmetric about $\frac{1}{2}$, implying
\[
C_{(a_n)} + C_{(b_n)} = C_{(a_n)} - C_{(b_n)} + 1,
\]
so studying sums of central Cantor sets is equivalent to studying their differences. Recent results on differences of Cantor sets from the perspective of achievement sets appear in \cite{Now23}, \cite{FN23}, \cite{Now23a}, \cite{Now25}. In particular, \cite{Now23} provides sophisticated sufficiency criteria for recognizing the topological type of sums of central Cantor sets and reproduces Pourbarat's results on sums of Cantor sets with constant dissection ratios \cite{Pour18}.

The classical Newhouse Gap Lemma \cite{New79}, \cite[Ch. 2]{Kr92} provides a sufficient condition for $C_1 + C_2$ to be an interval: $\tau(C_1) \cdot \tau(C_2) \ge 1$, where $\tau$ denotes thickness. Nowakowski observed that this condition is far from necessary \cite[Prop. 2.8]{Now23}.

\begin{theorem}
Given any $\epsilon > 0$, there exist sequences $(a_n), (b_n) \in (0,1)^\mathbb{N}$ such that $\tau(C_{(a_n)}) \cdot \tau(C_{(b_n)}) \le \epsilon$ and $C_{(a_n)} - C_{(b_n)} = [-1,1]$.
\end{theorem}

For central Cantor sets with constant ratios of dissection $\alpha, \beta \in (0,1)$, if $\frac{\ln \frac{1-\alpha}{2}}{\ln \frac{1-\beta}{2}}$ is irrational, then
$\textnormal{dim}_H\bigl(C_{(\alpha)}+C_{(\beta)}\bigr)\,=\,\max\bigl\{ \textnormal{dim}_HC_{(\alpha)}+\textnormal{dim}_HC_{(\beta)},\,1\,\bigr\}$ \cite{PS09}.
Using achievement set methods, Nowakowski showed that this equality does not hold in general for arbitrary Cantor sets \cite[Prop. 2.11]{Now23}. Additionally, his unpublished results provide simple sufficient conditions for the sum of two central Cantor sets to be a Cantor set of measure zero.

\begin{theorem}
Let $\sum x_n$ and $\sum y_n$ be fast convergent series of sum $1$, with remainders $r_n^x$ and $r_n^y$ respectively. If
\[
\lim_{n \to \infty} 3^n r_n^x = 0,
\]
then $E(x_n) + E(x_n)$ is a Cantor set of measure zero. If
\[
\lim_{n \to \infty} 3^n (r_n^x + r_n^y) = 0,
\]
then $E(x_n) + E(y_n)$ is a Cantor set of measure zero.
\end{theorem}

\begin{problem}
Can the constant $3$ in the first implication be replaced by a smaller number?
\end{problem}

\begin{problem}
For all central Cantor sets $C_{(a_n)}$, if $\text{int}(C_{(a_n)} - C_{(a_n)}) \neq \emptyset$, does it follow that $0 \in \text{int}(C_{(a_n)})$?
\end{problem}

This problem, if answered positively, would simplify criteria for ensuring that the difference of a central Cantor set is a Cantor set, as one would only need sequences of gaps tending to zero near $0$.

Another question is motivated by the observation of Moreira and Yoccoz that the Palis dichotomy generically holds for dynamically defined Cantor sets \cite{MY}.

\begin{problem}
Is the Palis dichotomy generically true for central Cantor sets?
\end{problem}

A recent arithmetic decomposition result \cite{MNP23} states:

\begin{theorem}
Every infinite achievable set is the algebraic sum of two achievable Cantor sets of Lebesgue measure zero.
\end{theorem}

This generalizes an earlier result for central Cantor sets \cite[Thm. 1]{PWT} and motivates further problems:

\begin{problem}
Can every infinite achievable set be decomposed into the sum of two achievable sets of Hausdorff dimension $0$?
\end{problem}

A positive answer seems plausible, in light of \cite[Prop. 2.11]{Now23}. Piotr Miska posed related questions about sums of achievement sets. The first asks whether there exists an achievable Cantor set $C$ of Hausdorff dimension $0$ such that $C + C$ is also a Cantor set. To formulate the answer, we introduce homogeneous Cantor sets and semi-fast convergent sequences.

Let $(n_k)$ be a sequence of integers with $n_k \ge 2$ and $(c_k)$ a sequence of positive numbers with $n_k c_k < 1$. Let $I_0$ be a closed interval $[a,b]$. Suppose $I_k$ consists of $i_k$ closed intervals $I_{k,i}$ of equal length. Define $I_{k+1}$ by removing $n_{k+1}-1$ open intervals from each $I_{k,i}$ so that the remaining $n_{k+1}$ closed intervals each have length $c_{k+1} |I_{k,i}|$. Then
\[
C_{(n_k),(c_k)} := \bigcap_{k=0}^\infty I_k = \bigcap_{k=0}^\infty \bigsqcup_{i=1}^{i_k} I_{k,i}
\]
is a \textsl{homogeneous Cantor set}. Central Cantor sets are a subfamily of homogeneous Cantor sets. Specifically, for $(a_n) \in (0,1)^\mathbb{N}$, we have $C_{(a_n)} = C_{(2), (\frac{1 - a_k}{2})}$. A homogeneous Cantor set $A$ with $\min A = 0$ corresponds to a semi-fast convergent sequence $(\alpha_k; N_k)$, i.e., $A = E(\alpha_k; N_k)$ \cite[pp. 1526-1527]{BFPW2}. For a homogeneous Cantor set $C_{(n_k),(c_k)}$ with fundamental interval $[0,\eta]$, we have $N_k = n_k - 1$ and $\alpha_k = \eta \prod_{i=1}^k c_i$. Its Hausdorff dimension is \cite[p. 804]{Lu}
\begin{equation}
\text{dim}_H C_{(n_k),(c_k)} = \liminf_{k \to \infty} \frac{\log(n_1 \cdots n_k)}{-\log(c_1 \cdots c_k)} = \liminf_{k \to \infty} \frac{\log((N_1+1) \cdots (N_k+1))}{-\log \alpha_k + \log \eta}.
\end{equation}

We can now answer Miska's first question positively. Consider $\alpha_k := \frac{1}{2} \left(\frac{1}{3^{(k-1)^2}} - \frac{1}{3^{k^2}}\right)$. Then $(\alpha_k)$ is fast convergent, so $E(\alpha_k)$ is a Cantor set. Similarly, $(\alpha_k; 2)$ is semi-fast convergent, so $E(\alpha_k; 2)$ is a Cantor set \cite[Thm. 16]{BFPW2}. Moreover,
\[
\text{dim}_H E(\alpha_k; 2) = \liminf_{k \to \infty} \frac{\log 3^k}{-\log \frac{1}{2} \left(\frac{1}{3^{(k-1)^2}} - \frac{1}{3^{k^2}}\right)} = 0.
\]
Since $E(\alpha_k) \subset E(\alpha_k; 2)$, we also have $\text{dim}_H E(\alpha_k) = 0$, providing a positive answer to Miska's first question. His second question remains open:

\begin{problem}
Let $E(x_n)$ be such that $E(x_n) + E(x_n)$ is a Cantor set. Is it true that
\[
\text{dim}_H(E(x_n)) \ge \frac{1}{2} \, \text{dim}_H(E(x_n) + E(x_n)) \,?
\]
\end{problem}

 \section{Cantovals, dual Cantorvals and their boundaries}

Given a compact set $A \subset \mathbb{R}$, we can fill in all the $A$-gaps and remove the interiors of all $A$-intervals. We denote the resulting set by $A^*$ and call it the \textsl{dual} of $A$. Formally,
\[
A^* = \overline{[\min A, \max A] \setminus A}.
\]
The dual of any Cantorval is a Cantorval. On the other hand, the boundary of any Cantorval is a Cantor set. It is known that the boundary of the Guthrie–Nymann Cantorval $E := E(3,2;\,\tfrac14)$ is an achievement set \cite[Example 7.2]{BGM18}, while its dual $E^*$ is not achievable \cite[Thm. 5.3]{BPW}.

Bartoszewicz posed the following question during the 5th Workshop on Postmodern Real Analysis:

\begin{problem}
Does there exist an achievement set $E$ such that both $\operatorname{Fr}\,E$ and $E^*$ are achievement sets?
\end{problem}

It is straightforward that the boundaries of achievable multi-interval sets need not be achievable. Of course, the boundaries of achievable Cantor sets are achievable, but the situation remains unclear for achievable Cantorvals.

\begin{problem}
Is the boundary of any achievable Cantorval an achievable set?
\end{problem}

These questions remain interesting even when restricted to multigeometric sequences. The next problem appears deceptively simple but highlights how little is known about achievable Cantorvals:

\begin{problem}
Do there exist two achievable Cantorvals $E(x_n)$ and $E(y_n)$ such that their union $E(x_n) \cup E(y_n)$ is a closed interval?
\end{problem}

Clearly, this can occur only if they share the same fundamental interval, that is, if $\sum x_n = \sum y_n$, under the standard assumption that all terms are positive.

Although the boundaries of Cantorvals, being Cantor sets, do not present topological difficulties, one can still ask about their size. The following problem has remained open for at least nine years:

\begin{problem}
Do there exist achievable Cantorvals with boundaries of positive Lebesgue measure?
\end{problem}

All known examples of achievable Cantorvals have boundaries of measure zero. This includes both multigeometric Cantorvals (the Guthrie–Nymann Cantorval \cite[Thm. 5.3]{BPW}, the Ferens Cantorvals \cite[Thm. 8]{BP}, the Guthrie–Nymann–Jones Cantorvals \cite[Thm. 3.1]{B19}) and non-multigeometric ones (the generalized Ferens Cantorvals \cite[Thm. 3.3]{MNP23}, the Marchwicki–Miska Cantorvals \cite[Thm. 9]{PWP}). Even the Kyiv Cantorvals \cite{VMPS19} share this property.

A Cantorval $E(x_n)$ is called \textsl{standard} if
\[
\varlimsup \frac{|P_k|}{r_k} > 0,
\]
where $|P_k|$ is the maximal length of $E_k$-intervals. Recently, Nowakowski and Prus-Wi\'{s}niowski \cite{NP} showed that all standard achievable Cantorvals have boundaries of Lebesgue measure zero.

\begin{problem}
Is every achievable Cantorval standard?
\end{problem}

A positive answer to this question would imply a negative answer to the previous one.

In another recent paper \cite{PK}, Karvatsky and Pratsiovytyi computed the Hausdorff dimension of the boundary of the Guthrie–Nymann Cantorval. They posed the following questions:

\begin{problem}
Is there an achievable Cantorval $E$ such that
\begin{itemize}
    \item[(i)] $\operatorname{dim}_H(\operatorname{Fr}\, E) = 1$?
    \item[(ii)] $\operatorname{dim}_H(\operatorname{Fr}\, E) = 0$?
\end{itemize}
\end{problem}

A crucial tool in proving the Guthrie–Nymann Classification Theorem is the Nymann–Saenz Lemma \cite[Lemma 2]{NS}.

\begin{problem}
Is there an analogue of the Nymann–Saenz Lemma for endpoints of interval components of an achievable Cantorval? More precisely, if $E(x_n)$ is a Cantorval and $a$ and $b$ are left endpoints of two distinct $E$-intervals, does there exist $\epsilon>0$ such that
\[
E \cap [a-\epsilon,a] \quad \text{is a translate of} \quad E \cap [b-\epsilon,b] \, ?
\]
\end{problem}

 \section{Achievement sets on the real plane}

 The concept of achievement sets extends naturally to any Banach space, in particular to \(\mathbb{R}^n\). Several earlier works have explored this direction (see \cite{BG15}, \cite{BGM18a}, \cite{GM18} and references therein). Fractal properties and dimensions of achievement sets were studied by Mor\'{a}n \cite{Mor89, Mor94}, and they can also be interpreted as ranges of purely atomic vector-valued measures \cite{LS08}. However, higher-dimensional achievement set theory is challenging and has not yet yielded spectacular results, mainly because there is no analog of the one-dimensional Guthrie-Nymann Classification Theorem. Consequently, a fundamental open problem in this subarea is the following:

\begin{problem}
Characterize achievement sets in \(\mathbb{R}^2\) up to homeomorphism.
\end{problem}

Proposition 3.1 in the recent paper \cite{KN24} lists the nine elementary topological types of achievement sets in the plane. Theorem 4.1 of the same paper states that, under mild assumptions on \((a_n)\), any set of \(P\)-sums \(S(P,(a_n))\) is a horizontal section of a suitably chosen two-dimensional achievement set \(E(x_n,y_n)\), which implies the existence of at least one additional topological type.

The following open problem is equivalent to Question 6.1 from \cite{KN24}:

\begin{problem}
Let \(E \subset \mathbb{R}^2\) be an achievable set. Is every vertical section \(E_x\) a set of \(P\)-sums?
\end{problem}

Multi-dimensional achievement sets appear naturally in modern number theory. A classical result of Erd\"{o}s and Straus states that the two-dimensional achievement set of \(\sum_{n=1}^\infty \left(\frac{1}{n}, \frac{1}{n+1}\right)\) has non-empty interior \cite[p.~65]{EG}. This was recently generalized to three dimensions by Kova\v{c} \cite{Kov25}. Moreover, some fundamental Kakeya observations on achievement sets were utilized in recent work by Kova\v{c} and Tao \cite{KT}. Thus, substantial development of achievement set theory may have direct implications for other important areas of mathematics.

\begin{problem}
Find reasonable sufficient conditions for \(\operatorname{int} E(x_n, y_n) \neq \emptyset\).
\end{problem}

Even the topological nature of the achievement sets of \(\sum_{n=1}^\infty (p^n, q^n)\) is not fully understood. Current observations are \cite[p.~221]{KN24}:
\begin{itemize}
    \item If at least one of the parameters \(p\) or \(q\) is less than or equal to \(\frac{1}{2}\), then \(E(p^n, q^n)\) is a Cantor set.
    \item If \(p = q \geq \frac{1}{2}\), then \(E(p^n, q^n)\) is a line segment.
\end{itemize}

\begin{problem}
Characterize \(E(p^n, q^n)\) for all \(p, q \in (0,1)\).
\end{problem}

It is straightforward that the Cartesian product of two one-dimensional achievement sets is itself a two-dimensional achievement set. By Theorem 4.1 in \cite{KN24}, however, there exist many two-dimensional achievement sets that are not Cartesian products of one-dimensional achievement sets, since there are numerous non-achievable sets of \(P\)-sums. For instance, every multi-interval subset of \([0,+\infty)\) is a set of \(P\)-sums, as recently proved by G\l\c{a}b and Marchwicki in \cite{GM26}.

\section{Cardinal function}

Let \(\mathrm{x} = (x_n)\) be an absolutely convergent series. Define \(f: \{0,1\}^\mathbb{N} \to \mathbb{R}\) by
\[
f((\varepsilon_n)) = \sum_{n=1}^\infty \varepsilon_n x_n.
\]
Then \(f\) is continuous, so \(f^{-1}(t)\) is closed and has cardinality in \(\{0,1,2,\dots,\omega,\mathfrak{c}\}\). The function \(F_\mathrm{x}(t) = \operatorname{card}(f^{-1}(t))\) is called the \emph{cardinality function} for \(\mathrm{x}\). The set
\[
U(E(x_n)) := \{t \in E(x_n) : F_\mathrm{x}(t) = 1\}
\]
is called the \emph{set of uniqueness} of \(E(x_n)\) and is a \(\mathcal{G}_\delta\) set \cite{NP}. Note that \(\{0, \sum_{n=1}^\infty x_n\} \subseteq U(E(x_n))\), and equality occurs if and only if \(E(x_n)\) is the interval \([0, \sum_{n=1}^\infty x_n]\).

In \cite{GM20}, the following relations between Kakeya conditions and cardinal functions were established:
\begin{itemize}
    \item If \(x_n > r_n\) for every \(n\), then \(U(E(x_n)) = E(x_n)\).
    \item If \(x_n \le r_n\) for every \(n\) and \(x_n < r_n\) for only finitely many \(n\), then \(F_\mathrm{x}\) is bounded.
    \item If \(x_n \le r_n\) for every \(n\) and \(x_n < r_n\) for infinitely many \(n\), then \(F_\mathrm{x}\) is unbounded, and there exists \(t\) with \(F_\mathrm{x}(t) \in \{\omega, \mathfrak{c}\}\).
\end{itemize}

In \cite{GM23}, it was shown that \(U(GN)\) is residual in \(GN\), which implies that \(GN \setminus U(GN)\) is topologically small. However, \(U(E(x_n))\) never contains an interval.

\begin{problem}
For which sets \(C \subseteq \{0,1,2,\dots,\omega,\mathfrak{c}\}\) does there exist a sequence \(\mathrm{x}\) such that the cardinal function \(F_\mathrm{x}\) has range exactly \(C\)?
\end{problem}

Partial answers appear in \cite{MZ}. The following specific question remains open:

\begin{problem}[\cite{GM20}]
Does there exist a sequence \(\mathrm{x}\) such that the range of \(F_\mathrm{x}\) contains an unbounded sequence of natural numbers but neither \(\omega\) nor \(\mathfrak{c}\)?
\end{problem}

If the set \(\{t \in \mathbb{R} : F_\mathrm{x}(t) = \omega \text{ or } F_\mathrm{x}(t) = \mathfrak{c}\}\) is non-empty, then it is dense in \(E(x_n)\) and therefore infinite \cite{GM23}.

\begin{problem}[\cite{GM20}]
Is it true that either \(\{t \in \mathbb{R} : F_\mathrm{x}(t) = \mathfrak{c}\}\) (or \(\{t \in \mathbb{R} : F_\mathrm{x}(t) = \omega\}\)) is empty or contains a Cantor set?
\end{problem}

A sequence \((x_n)\) is a \emph{minimal representation} of a Cantorval or a multi-interval set if removing any infinite subsequence changes the topological type of the achievement set.

\begin{problem}
It is known that if \(E\) is an achievable multi-interval set, then \((x_n)\) is a minimal representation for \(E\) if and only if \(x_n < r_{n+1}\) for only finitely many \(n\) \cite[Thm. 7.8]{GM23}. Does this characterization remain valid when \(E\) is a Cantorval?
\end{problem}

\section{Center of distances}

The \emph{center of distances} is a metric invariant introduced by Bielas, Plewik, and Walczy\'{n}ska \cite{BPW} as a tool to detect non-achievability. Given a metric space \((X,d)\), the center of distances of \(A \subset X\) is
\[
\mathcal{C}(A) := \{\alpha \ge 0 : \forall x \in A, \exists y \in A \text{ with } d(x,y) = \alpha\}.
\]
Clearly, \(0 \in \mathcal{C}(A) \subset [0,+\infty)\) for any \(A \subset \mathbb{R}\). The recent paper \cite{BFHLP} studies basic properties of \(\mathcal{C}\) as a map from the hyperspace of compact subsets of \([0,1]\) into itself. In particular, points of continuity of this map are exactly those compact \(A \subset [0,1]\) with \(\mathcal{C}(A) = \{0\}\) \cite[Thm. 3.7]{BFHLP}. Another striking result is the following \cite[Cor. 3.9]{BFHLP}:

\begin{theorem}
The family of achievement sets is of the first category in the hyperspace of compact subsets.
\end{theorem}
Most recently, it has been shown that the set of achievement sets is closed in the hyperspace and hence nowhere dense \cite[Thm. 4.5]{NPT}.

While the center of distances is practically useless in higher dimensions, it can be replaced by the promising concept of the \emph{spectre of a set}, recently introduced by Kula and Nowakowski \cite{KN24}.

Another important new result is the following \cite[Thm. 3]{K25}, which answers positively a question posed by M. Filipczak one year earlier:

\begin{theorem}
Every set \(B \subset \mathbb{R}\) with \(0 \in B \subset [0,+\infty)\) is the center of distances of some set \(A \subset \mathbb{R}\).
\end{theorem}

The proof uses transfinite induction, so it does not provide information on descriptive properties of \(A\) for a given \(B\). For example:

\begin{problem}
Is it true that for any Borel set \(0 \in A \subset [0,\infty)\), there exists a Borel set \(B \subseteq \mathbb{R}\) with \(\mathcal{C}(B) = A\)?
\end{problem}

Of particular interest is the following:

\begin{problem}
Assume \(0 \in C \subseteq [0,\infty)\) is compact. Does there exist a compact \(K \subseteq \mathbb{R}\) with \(\mathcal{C}(K) = C\)?
\end{problem}

This problem was discussed during the 1st Workshop on Postmodern Real Analysis in 2020, but only an unsatisfactory partial answer was obtained \cite[Thm. 4.10]{BFHLP}.

Given a central Cantor set \(C = C_{(a_n)}\), the union \(G_n\) of all intervals removed at the \(n\)-th step is well-defined. Define
\[
\widehat{C} := C \cup \bigcup_{n \in 2\mathbb{N}-1} G_n,
\]
which is a Cantorval called the \emph{adjoint Cantorval} of \(C\). For example, \(\widehat{C_{(1/3)}}\) is the model Cantorval. One elementary but surprisingly resilient problem is:

\begin{problem}
Compute \(\mathcal{C}(\widehat{E(q^n)})\) for \(q \in (0, \frac12)\).
\end{problem}

 It is almost certain that \(\mathcal{C}(\widehat{E(q^n)}) = \{0\}\) for large \(q\), and it is known that \(\mathcal{C}(\widehat{E(q^n)})\) contains intervals for \(q < \frac14\).\footnote{Readers interested in this problem are encouraged to contact Jacek Marchwicki at \texttt{jacek.marchwicki@uwm.edu.pl}.}

Since any central Cantor set \(C\) is the achievement set of a unique fast convergent series \(\sum x_n\), it follows that
$
\{0\} \cup \{x_n : n \in \mathbb{N}\} \subset \mathcal{C}(C)
$ \cite[Prop. 3.1]{BPW}.
The recent paper \cite{BBFP} investigates when a central Cantor set has minimal center of distances, i.e., when
\[
\{0\} \cup \{x_n : n \in \mathbb{N}\} = \mathcal{C}(C),
\]
and Banakiewicz completed this work in \cite[Thm. 2]{B23}:

\begin{theorem}
     Let $\sum x_n$ be a fast convergent series. Then the center of distances of the central Cantor set $E(x_n)$ is not minimal if and only if at least one of the three following conditions holds
     \begin{itemize}
         \item[\textnormal{(B1)}]\ $\exists\ n\in\mathbb N\quad x_{n-2}=4x_n,\ x_{n-1}=2x_n;$
         \item[\textnormal{(B2)}]\ $\exists\ n\in\mathbb N\quad x_{n-3}=9x_n, x_{n-2}=5x_n,\ x_{n-1}=2x_n;$
          \item[\textnormal{(B3)}]\ $\exists\ n\in\mathbb N\quad x_{n-3}=10x_n, x_{n-2}=6x_n,\ x_{n-1}=2x_n.$
     \end{itemize}
 \end{theorem}

The last problem presented here concerns adjoint Cantorvals and was posed by Kula during the 5th Workshop on Postmodern Real Analysis in 2023:

\begin{problem}
Do there exist fast convergent sequences \((x_n)\) such that the adjoint Cantorval \(\widehat{E(x_n)}\) is achievable?
\end{problem}

We conjecture that the answer is negative, because it is so in the case of fast convergent geometric sequences as it was shown in \cite{BFGPWS}.


\begin{thebibliography}{99999}
\bibitem{BBFS}
T. Banakh, A. Bartoszewicz, M. Filipczak, E. Szymonik, \textit{Topological and measure properties of some self-similar sets}, Topol. Methods Nonlinear Anal., 46(2) (2015), 1013--1028
\bibitem{B19}
M. Banakiewicz, \textit{The Lebesgue measure of some $M$-Cantorval}, J. Math. Anal. Appl. 471(2019)  170-179
\bibitem{B23}
M. Banakiewicz, \textit{The center of distances of central Cantor sets}, Results Math. (2023) 78:234  https://doi.org/10.1007/s00025-023-02012-3
\bibitem{BBFP}
M. Banakiewicz, A. Bartoszewicz, M. Filipczak, F. Prus-Wi\'{s}niowski, \textit{Center of distances and central Cantor sets}, Results Math. (2022) 77:196  https://doi.org/10.1007/s00025-022-01725-1
\bibitem{BP}
M. Banakiewicz, F. Prus-Wi\'{s}niowski, \textit{M-Cantorvals of Ferens type}, Math. Slovaca, 67(4)(2017), 1-12
\bibitem{BFGPWS} A. Bartoszewicz, M. Filipczak, S. G\l\c{a}b, F. Prus-Wi\'sniowski, J. Swaczyna, \textit{On generating regular Cantorvals connected with geometric Cantor sets}, Chaos, Solitons and Fractals 114 (2018) 468--473
\bibitem{BFHLP}
A. Bartoszewicz, M. Filipczak, G. Horbaczewska, S. Lindner, F. Prus-Wi\'{s}niowski, \textit{On the operator of center of distances between the spaces of closed subsets of the real line}, Topol. Metods Nonlinear Anal.  63(2) (2024) 413–427
DOI: 10.12775/TMNA.2023.023
\bibitem{BFPW1} A. Bartoszewicz, M. Filipczak, F. Prus-Wi\'{s}niowski, \textit{Topological and algebraic aspects of subsums of series}, Traditional and present-day topics in real analysis, 345--366, Faculty of Mathematics and Computer Science. University of \L \'od\'z, \L \'od\'z, 2013.
\bibitem{BFPW2} A. Bartoszewicz, M. Filipczak, F. Prus-Wi\'{s}niowski, \textit{Semi-fast convergent sequences and $k$-sums of central Cantor sets}, European J. Math. 6(2020), 1523-1536
\bibitem{BFS}
A. Bartoszewicz, M. Filipczak, E. Szymonik, \textit{Multigeometric sequences and Cantorvals}, Cent. Eur. J. Math. 12(7) (2014), 1000--1007.
\bibitem{BG15} A. Bartoszewicz, S. G\l\c{a}b. \emph{Achievement sets on the plane — Perturbations of geometric and multigeometric series}. Chaos Solit. Fractals \textbf{77} (2015), 84--93.
\bibitem{BGM18}
A. Bartoszewicz, S. G\l \c{a}b, J. Marchwicki, \textit{Recovering a purely atomic finite measure from its range}, J. Math. Anal. Appl. 467 (2018) 825-841.
\bibitem{BGM18a}
A. Bartoszewicz, S. G\l\c{a}b, J. Marchwicki. \emph{Achievement sets of conditionally convergent series}. Colloq. Math. \textbf{152} (2018), 235--254.
\bibitem{BPW}  W. Bielas, S. Plewik, M. Walczy\'{n}ska, \textit{On the center of distances}, Eur. J. Math. 4 (2018), no. 2, 687--698.
\bibitem{BKP} Bielas, W. Kula, M. Plewik, Sz. \textit{On compact subsets of the reals}, Topology Appl.346(2024), Paper No. 108854, 10 pp.
\bibitem{CJPPS}M. C\"ornyei, T. Jordan, M. Pollicott, D. Preiss, B. Solomyak, \textit{Positive-measure self-similar sets without interior}, Ergodic Theory
Dynam. Systems. 26:3 (2006), 755--758.
\bibitem{EG}
P. Erd\"{o}s, R.L. Graham, \textit{Old and new problems and results in combinatorial number theory},
volume 28 of Monographies de L’Enseignement Math´ematique. Universit´e de Gen`eve, L’Enseignement
Math´ematique, Geneva, 1980.
\bibitem{F}
C. Ferens, \textit{On the range of purely atomic probability
measures}, Studia Math., \textbf{77(3)} (1984), 261--263.
\bibitem{FN23}
T. Filipczak, P. Nowakowski, \textit{Conditions for the difference set of a
central Cantors set to be a Cantorval}, Results Math 78, 166 (2023). https://doi.org/10.1007/s00025-023-01940-4
\bibitem{GK25} S. G\l\c{a}b, M. Kula, \textit{Kenyon Theorem Revisited}, arXiv:2508.19454v1 [math.GN] 26Aug2025
\bibitem{GM18}
S. G\l\c{a}b, J. Marchwicki, \textit{Levy-Steinitz theorem and achievement sets of conditionally
convergent series on the real plane}, J. Math. Anal. Appl., 459(1) (2018) 476–489
doi:10.1016/j.jmaa.2017.10.034
\bibitem{GM20} S. G\l\c{a}b, J. Marchwicki, \textit{Cardinal functions of purely atomic measures}, Results in Math. 75 (4), art. nr.141 (2020), 26 pp.
\bibitem{GM23} S. G\l\c{a}b, J. Marchwicki, \textit{Set of uniqueness for Cantorvals}, Results in Math. 78(9) (2023) DOI:10.1007/s00025-022-01777-3
\bibitem{GM26}
S. G\l\c{a}b, J. Marchwicki, \textit{On arithmetic sums of Cantor sets: \(P\)-sums vs. achievement sets}, accepted in the Rivista di Matematica della Universit\`{a} di Parma
\bibitem{GN88} J.A. Guthrie, J.E. Nymann, \textit{The topological structure of the set of subsums of an infinite series}, Colloq. Math. 55:2 (1988), 323--327.
\bibitem{Ho41}
Hornich, H., \textit{\"{U}ber beliebige Teilsummen absolut konvergenter Reihen}, Monatsh. Math. Phys. 49(1941) 316-320
\bibitem{Jones} R. Jones, \textit{Achievement sets of sequences}, Am. Math. Mon. 118:6 (2011), 508--521.
\bibitem{Kakeya} S. Kakeya, \textit{On the partial sums of an infinite series}, T\^{o}hoku c. Rep. 3 (1914), 159--164.
\bibitem{Kakeya2} S. Kakeya, \textit{On the set of partial sums of an infinite series}, Proc. Tokyo Math.-Phys. Soc., 2nd series, 7(1914), 250-251,  \text{\path{ https://doi.org/10.11429/ptmps1907.7.14_250}}
\bibitem{Kar24} D. Karvatskyi, A. Murillo, A.  Viruel, \textit{The Achievement Set of Generalized Multigeometric Sequences.} Results Math 79, 132 (2024). https://doi.org/10.1007/s00025-024-02158-8

\bibitem{Kov25}
Kovač, V. (2025), \textit{On the Set of Points Represented by Harmonic Subseries}, Amer. Math. Monthly, 1–17. https://doi.org/10.1080/00029890.2025.2540753
\bibitem{KT}
V. Kova\v{c}, T. Tao, \textit{On several irrationality problems for Ahmes series}, Acta Math. Hungar. 175(2) (2025), 572-608. https://doi.org/10.1007/s10474-025-01528-0
\bibitem{Kr92}
R. Kraft, \textit{Intersections of thick Cantor sets}, Memoirs of the Amer. Math. Soc., vol. 97(468), Providence, Rhode Islanc, 1992
\bibitem{K25}
M. Kula, \textit{Center of distances and Bernstein sets}, Real Anal. Exchange 50(1) (2025) 207-212,  DOI: 10.14321/realanalexch.1739330962
\bibitem{KN24}
M. Kula, P. Nowakowski, \textit{Achievement sets of series in $\mathbb R^2$}, Results Math. (2024) 79:221 https://doi.org/10.1007/s00025-024-02239-8
\bibitem{LS08}
Laltanpuia,  A. I.  Singh, \textit{ Preservation of the range of a vector measure under shortenings of the
domain}, J. Convex Anal., 15(2)(2008), 313–324.
\bibitem{Lu}
S. Lu, \textit{The Housdorff dimension and measure of some Cantor sets}, Real Anal. Exchange 25(2) (1999/2000) 799-808
\bibitem{M24}
M. Moroz, \textit{A counterexample to the Karvatskyi-Pratsiowytyi conjecture concerning the achievement set of an intermediate series}, arXiv:2412.00042v1 [mathGM]
22Nov2024
\bibitem{MM21}
J. Marchwicki, J. Miska, \textit{On Kakeya conditions for achievement sets}, Result in Math. (2021) 76:181; //doi.org./10.1007/s00025-021-01479-2
\bibitem{MNP23}
J. Marchwicki, P. Nowakowski, F. Prus-Wiśniowski, \textit{Algebraic sums of achivable sets involving Cantorvals}, arXiv:2309.01589v1 [math.CA]
\bibitem{MZ} J. Marchwicki, B. \.Zmija, \textit{Bounded ranges of cardinal functions}, arXiv:2508.13016 [math.CA]
\bibitem{MO} P. Mendes, F. Oliveira, \textit{On the topological structure of the arithmetic sum of two Cantor sets}, Nonlinearity. 7 (1994), 329--343
\bibitem{MPWP}
        Miska, P., Prus-Wiśniowski, F. \& Ptak, J. \textit{More on Kakeya conditions for achievement sets}, Results Math 78, 113 (2023). https://doi.org/10.1007/s00025-023-01890-x
\bibitem{Mor89}
M. Mor\'{a}n, \textit{Fractal series}, Mathematika, 36(2) (1989) 334–348, doi:10.1112/S0025579300013176
\bibitem{Mor94}
M. Mor\'{a}n, \textit{Dimension functions for fractal sets associated to series}, Proc. Amer. Math. Soc.,
120(3) (1994) 749–754,  doi:10.2307/2160466
 \bibitem{MY}
 C. Moreira, J. Yoccoz, \textit{Stable intersections of Cantor sets with large Hausdorff dimension}, Ann. of Math. 154(2001), 45 - 96
 \bibitem{New79}
 S. Newhouse, \textit{The abundance of wild hyperbolic sets and nonsmooth stable sets for diffeomorphisms}, Inst. Hautes \'{E}tudes Sci. Publ. Math. 50 (1979), 101-105
\bibitem{Nit11} Z. Nitecki, \textit{The subsum set of a null sequence}, arXiv:1106.3779v1 [math.HO] 19Jun2011
\bibitem{Nit13} Z. Nitecki, \textit{Subsum sets: intervals, Cantor sets, and Cantorvals}, arXiv:1106.3779v2 [math.HO] 8Jul2013
\bibitem{Nit15} Z. Nitecki, \textit{Cantorvals and subsum sets of null sequences}, Amer. Math. Monthly 122 (2015), 862--870.
\bibitem{Now23}
P. Nowakowski, \textit{When the algebraic difference of two central Cantor sets is an interval ?}, Ann. Fennici Math. 48(2023), 163-185
\bibitem{Now23a}
P. Nowakowski, \textit{Characterization of the algebraic difference of special affine Cantor sets }, Topol. Methods in Nonlinear Anal. 64(1) (2024)  295 - 316   	https://doi.org/10.12775/TMNA.2023.057
\bibitem{Now25}
P. Nowakowski, \textit{Conditions for the difference set of a central Cantor set to be a Cantorval. Part II}, Indag. Math. 36 (2025) 1223–1244
\bibitem{NP}
P. Nowakowski, F. Prus-Wi\'{s}niowski, \textit{The Lebesgue measure of boundaries of multigeometric Cantorvals}, arXiv: 2510.21878v1 [mathDS] 23 Oct 2025
\bibitem{NPT}
P. Nowakowski, F. Prus-Wi\'{s}niowski, F. Turobo\'{s}, \textit{Spectre operator, achievement sets and sets of P-sums in a hyperspace of compact sets}  	arXiv:2512.11803 [math.GN] 2 Nov 2025
\bibitem{NS}
J.E. Nymann, R.A. S\'{a}enz, \textit{On a paper of Guthrie and Nymann on subsums of infinite series}, Colloq. Math. 83(2000) 1-4
\bibitem{PS09}
Y. Peres, P. Shmerkin, \textit{Resonanse between Cantor sets}, Ergod. Theory Dyn. Syst. 29 (2009), 201-22
\bibitem{Pour18}
M. Pourbarat, \textit{On the arithmetic difference of middle Cantor sets}, Discrete Contin. Dyn. Syst. 38(9) (2018), 4259-4278
\bibitem{PK23}
M. Pratsiovytyi, D. Karvatskyi, \textit{Cantorvals as sets of subsums for a series
related with trigonometric functions}, Proc.
Intern. Geom. Center
Vol. 16, no. 3-4 (2023) pp. 262–271
\bibitem{PK}
M. Pratsiovytyi, D. Karvatskyi, \textit{Fractal analysis of Guthrie-Nymann's set and its generalisations},     DOI: 10.48550/arXiv.2405.16576
\bibitem{PR25} M. Pratsiovytyi, S. Ratushniak, \textit{Properties of a disturbed binary series}. Bukovinian Math. J. 13(1) (2025), 109-117.
\bibitem{PWP}
F. Prus-Wi\'{s}niowski, J. Ptak,\textit{Achievable Cantorvals almost without reversed Kakeya conditions}, arXiv:2412.08768v1 [math.CA] 11 Dec 2024
\bibitem{PWT}
F. Prus-Wi\'{s}niowski, F. Tulone, \textit{The arithmetic decomposition of central Cantor sets}, J. Math. Anal. Appl. 467(2018) 26-31
\bibitem{W1}
D. Waterman, \textit{$\Lambda$-bounded variation: recent results and unsolved problems}, Real Anal. Exchange 4(1978-79) 69-75
\bibitem{W2}
D. Waterman, \textit{Generalized bounded variation: recent results and open questions}, Real Anal. Exchange 5(1979-80) 148-150
\bibitem{WS}
A.D. Weinstein, B.E. Shapiro, \textit{On the structure of a set of $\overline{\alpha}$-representable numbers}, Izv. Vys\v{s}. U\v{c}ebn. Zaved. Matematika. 24 (1980), 8--11.
\bibitem{VMPS19}
Vinishin Y., Markitan V., Pratsiovytyi M., Savchenko I., \textit{Positive series, whose sets of subsums are Cantorvals}, Proc. International Geom. Center. 2019, 12 (2), 26-42. (in Ukrainian);
https://doi.org/10.15673/tmgc.v12i2.1455

\end{thebibliography}
\end{document}